\documentclass[11pt,leqno]{article}

\usepackage[active]{srcltx}

\usepackage{amsfonts,latexsym,amsmath,amssymb,amsthm}
\usepackage{fullpage}
\newtheorem{theorem}{Theorem}[section]
\newtheorem{lemma}[theorem]{Lemma}

\newtheorem{question}[theorem]{Question}
\newtheorem{proposition}[theorem]{Proposition}

\newtheorem{remark}[theorem]{Remark}

\theoremstyle{definition}

\theoremstyle{remark}

\newtheorem*{note*}{Note}

\numberwithin{equation}{section}

\makeatletter
\newcommand{\rank}{\mathop{\operator@font rank}}
\newcommand{\conv}{\mathop{\operator@font conv}}
\newcommand{\vol}{\mathop{\operator@font vol}}
\newcommand{\onetagright}{\tagsleft@false}
\makeatother

\newcommand{\ls}{\leqslant}
\newcommand{\gr}{\geqslant}
\renewcommand{\epsilon}{\varepsilon}
\newcommand{\prend}{$\quad \hfill \Box$}
\usepackage{color}

\usepackage{hyperref}

\begin{document}
\small

\title{\bf Geometry of random sections of isotropic convex bodies}

\medskip

\author{Apostolos Giannopoulos, Labrini Hioni and Antonis Tsolomitis}

\date{}

\maketitle

\begin{abstract}
\footnotesize Let $K$ be an isotropic symmetric convex body in ${\mathbb
R}^n$. We show that a subspace $F\in G_{n,n-k}$ of codimension $k=\gamma
n$, where $\gamma\in (1/\sqrt{n},1)$, satisfies $$K\cap F\subseteq
\frac{c}{\gamma }\sqrt{n}L_K (B_2^n\cap F)$$ with probability greater than
$1-\exp (-\sqrt{n})$. Using a different method we study the same question
for the $L_q$-centroid bodies $Z_q(\mu )$ of an isotropic log-concave
probability measure $\mu $ on ${\mathbb R}^n$. For every $1\ls q\ls n$ and
$\gamma\in (0,1)$ we show that a random subspace $F\in G_{n,(1-\gamma )n}$
satisfies $Z_q(\mu )\cap F\subseteq c_2(\gamma )\sqrt{q}\,B_2^n\cap F$. We
also give bounds on the diameter of random projections of $Z_q(\mu )$ and
using them we deduce that if $K$ is an isotropic convex body in ${\mathbb R}^n$
then for a random subspace $F$ of dimension $(\log n)^4$ one has that
all directions in $F$ are sub-Gaussian with constant $O(\log^2n)$.
\end{abstract}

\section{Introduction}

A convex body $K$ in ${\mathbb R}^n$ is called isotropic if it has volume $|K|=1$, its center of mass is at the origin
(we call these convex bodies ``centered"), and its inertia matrix is a multiple of the identity matrix: there exists a constant $L_K >0$ such that
\begin{equation}\label{eq:intro-1}\int_K\langle x,\theta\rangle^2dx =L_K^2\end{equation}
for every $\theta $ in the Euclidean unit sphere $S^{n-1}$. For every centered convex body $K$ in ${\mathbb R}^n$
there exists an invertible linear transformation $T\in GL(n)$ such that $T(K)$ is isotropic. This isotropic image of $K$ is
uniquely determined up to orthogonal transformations. A well-known problem in asymptotic convex geometry asks if there exists an absolute constant $C_1>0$ such that
\begin{equation}\label{eq:intro-2}L_n:= \max\{ L_K:K\ \hbox{is isotropic in}\ {\mathbb R}^n\}\ls C_1\end{equation}
for all $n\gr 1$ (see Section 2 for background information on isotropic convex bodies and log-concave probability
measures). Bourgain proved in \cite{Bourgain-1991} that $L_n\ls c\sqrt[4]{n}\log\! n$, and Klartag \cite{Klartag-2006}
improved this bound to $L_n\ls c\sqrt[4]{n}$. A second proof of Klartag's bound appears in \cite{Klartag-EMilman-2012}.

Recall that the inradius $r(K)$ of a convex body $K$ in ${\mathbb R}^n$ with $0\in {\rm int}(K)$ is the largest $r>0$ for which $rB_2^n\subseteq K$, while
the radius $R(K):=\max\{\| x\|_2:x\in K\}$ of $K$ is the smallest $R>0$ for
which $K\subseteq RB_2^n$. It is not hard to see that the inradius
and the radius of an isotropic convex body $K$ in ${\mathbb R}^n$ satisfy the bounds $c_1L_K\ls r(K)\ls R(K)\ls c_2nL_K$,
where $c_1,c_2>0$ are absolute constants. In fact, Kannan, Lov\'{a}sz and Simonovits \cite{Kannan-Lovasz-Simonovits-1995} have proved
that
\begin{equation}\label{eq:intro-3}R(K)\ls (n+1)L_K.\end{equation}

\smallskip

\noindent {\bf Radius of random sections of isotropic convex bodies.} The first question that we discuss in this article
is to give sharp upper bounds for the radius of a random $(n-k)$-dimensional section of $K$. A natural ``guess" is that the following
question has an affirmative answer.

\begin{question}\label{question}There exists an absolute constant $\overline{c}_0>0$ with the following property: for every isotropic convex body $K$ in ${\mathbb R}^n$
and for every $1\ls k\ls n-1$, a random subspace $F\in G_{n,n-k}$ satisfies
\begin{equation}\label{eq:question}R(K\cap F)\ls \overline{c}_0\sqrt{n/k}\,\sqrt{n}L_K.\end{equation}
\end{question}

It was proved in \cite{Litvak-Milman-Pajor-1999} that if $K$ is a symmetric convex body in ${\mathbb R}^n$ then a random
$F\in G_{n,n-k}$ satisfies
\begin{equation}\label{eq:intro-111}R(K\cap F)\ls c(n/k)^{3/2}\tilde{M}(K),\end{equation}
where $c>0$ is an absolute constant and
\begin{equation}\tilde{M}(K):=\frac{1}{|K|}\int_K\|x\|_2dx.\end{equation}
In the case of an isotropic convex body one has $|K|=1$ and
\begin{equation}\tilde{M}(K)\ls \left (\int_K\|x\|_2^2dx\right )^{1/2}=\sqrt{n}L_K,\end{equation}
therefore \eqref{eq:intro-111} implies that a random $F\in G_{n,n-k}$ satisfies
\begin{equation}\label{eq:intro-112}R(K\cap F)\ls \overline{c}_1(n/k)^{3/2}\sqrt{n}L_K,\end{equation}
where $\overline{c}_1>0$ is an absolute constant.

Our first main result shows that one can have a bound of the order of $\gamma^{-1}\sqrt{n}L_K$ when the codimension $k$ is greater than $\gamma n$.

\begin{theorem}\label{th:intro-2}Let $K$ be an isotropic symmetric convex body in ${\mathbb R}^n$ and let $1\ls k\ls n-1$.
A random subspace $F\in G_{n,n-k}$ satisfies
\begin{equation}\label{eq:intro-5}R(K\cap F)\ls \frac{\overline{c}_0n}{\max\{k,\sqrt{n}\}}\sqrt{n}L_K\end{equation}
with probability greater than $1-\exp (-\sqrt{n})$, where $\overline{c}_0>0$ is an absolute
constant.
\end{theorem}

The proof is given in Section 3. Note that Theorem \ref{th:intro-2} gives non-trivial information when $k>\sqrt{n}$. In this case, writing $k=\gamma n$
for some $\gamma\in (1/\sqrt{n},1)$ we see that
\begin{equation}\label{eq:intro-55}R(K\cap F)\ls \frac{\overline{c}_0}{\gamma }\sqrt{n}L_K\end{equation}
with probability greater than $1-\exp (-\sqrt{n})$ on $G_{n,(1-\gamma )n}$. The result of \cite{Litvak-Milman-Pajor-1999}
establishes a $\gamma^{-3/2}$-dependence on $\gamma =k/n$.

A standard approach to Question \ref{question} would have been to combine the low $M^{\ast }$-estimate with an upper bound for the mean width
\begin{equation}\label{eq:intro-7}w(K):=\int_{S^{n-1}}h_K(x)\,d\sigma (x),\end{equation}
of an isotropic convex body $K$ in ${\mathbb R}^n$, that is, the $L_1$-norm of the support function of $K$ with respect to the Haar
measure on the sphere. This last problem was open for a number of years. The upper bound $w(K)\ls cn^{3/4}L_K$ appeared in the Ph.D.
Thesis of Hartzoulaki \cite{Hartzoulaki-thesis}. Other approaches leading to the same bound can be found in Pivovarov \cite{Pivovarov-2010a}
and in Giannopoulos, Paouris and Valettas \cite{Giannopoulos-Paouris-Valettas-2012b}. Recently, E.~Milman showed in \cite{EMilman-2014} that
if $K$ is an isotropic symmetric convex body in ${\mathbb R}^n$ then
\begin{equation}\label{eq:intro-8}w(K)\ls c_3\sqrt{n}(\log n)^2L_K.\end{equation}
In fact, it is not hard to see that his argument can be generalized to give the same estimate in the
not necessarily symmetric case. The dependence on $n$ is optimal up to the logarithmic term. From the sharp version of V.~Milman's low $M^{\ast }$-estimate
(due to Pajor and Tomczak-Jaegermann \cite{Pajor-Tomczak-1986}; see \cite[Chapter 7]{AGA-book}
for complete references) one has that, for every $1\ls k\ls n-1$, a subspace $F\in G_{n,n-k}$ satisfies
\begin{equation}\label{eq:intro-9}R(K\cap F)\ls c_4\sqrt{n/k}\,w(K)\end{equation}
with probability greater than $1-\exp (-c_5k)$, where $c_4, c_5>0$ are absolute
constants. Combining \eqref{eq:intro-9} with E.~Milman's theorem we obtain the folowing estimate:
\begin{quote}{\sl Let $K$ be an isotropic symmetric convex body in ${\mathbb R}^n$. For every $1\ls k\ls n-1$, a subspace
$F\in G_{n,n-k}$ satisfies
\begin{equation}\label{eq:emanuel-bound}R(K\cap F)\ls \frac{\overline{c}_2n(\log n)^2L_K}{\sqrt{k}}\end{equation}
with probability greater than $1-\exp (-\overline{c}_3k)$, where $\overline{c}_2, \overline{c}_3>0$ are absolute
constants.}
\end{quote}
\noindent Note that the upper bound of Theorem \ref{th:intro-2} has some advantages when compared to \eqref{eq:emanuel-bound}:
If $k$ is proportional to $n$ (say $k\gr\gamma n$ for some $\gamma\in (1/\sqrt{n},1)$) then Theorem \ref{th:intro-2} guarantees that
$R(K\cap F)\ls c(\gamma )\sqrt{n}L_K$ for a random $F\in G_{n,n-k}$. More generally, for all $k\gr \frac{c_6n}{(\log n)^4}$ we have
\begin{equation}\label{eq:intro-10}\frac{\overline{c}_0n\sqrt{n}}{\max\{k,\sqrt{n}\}}\ls \frac{\overline{c}_2n(\log n)^2}{\sqrt{k}},\end{equation}
and hence the estimate of Theorem \ref{th:intro-2} is stronger than \eqref{eq:emanuel-bound}. Nevertheless, we emphasize that our
bound is not optimal and it would be very interesting to decide whether \eqref{eq:question} holds true; this would be optimal for all $1\ls k\ls n$.

\medskip

\noindent {\bf Radius of random sections of $L_q$-centroid bodies and their polars.}
In Section 4 we study the diameter of random sections of the $L_q$-centroid bodies $Z_q(\mu )$ of an isotropic log-concave probability measure
$\mu $ on ${\mathbb R}^n$. Recall that a measure $\mu$ on $\mathbb R^n$ is called log-concave if $\mu(\lambda A+(1-\lambda)B)
\gr \mu(A)^{\lambda}\mu(B)^{1-\lambda}$ for any compact subsets $A$ and $B$ of ${\mathbb R}^n$ and any $\lambda \in (0,1)$. A function
$f:\mathbb R^n \rightarrow [0,\infty)$ is called log-concave if its support $\{f>0\}$ is a convex set and the restriction of $\log{f}$ on it is concave.
It is known that if a probability measure $\mu $ is log-concave and $\mu (H)<1$ for every hyperplane $H$, then $\mu $ is absolutely
continuous with respect to the Lebesgue measure and its density
$f_{\mu}$ is log-concave; see \cite{Borell-1974}. Note that if $K$ is a convex body in $\mathbb R^n$ then the Brunn-Minkowski inequality implies that the indicator function
$\mathbf{1}_{K} $ of $K$ is the density of a log-concave measure.

We say that a log-concave probability measure $\mu $ on ${\mathbb R}^n$
is isotropic if its barycenter $\textrm{bar}(\mu )$ is at the origin and
\begin{equation*}\int_{{\mathbb R}^n}\langle x,\theta\rangle^2\,d\mu (x)=1\end{equation*}
for all $\theta\in S^{n-1}$. Note that the normalization is different from the one
in \eqref{eq:intro-1}; in particular, a centered convex body $K$ of volume $1$ in ${\mathbb R}^n$ is isotropic 
if and only if the log-concave probability measure $\mu_K$ with density
$x\mapsto L_K^n\mathbf{1}_{K/L_K}(x)$ is isotropic.

The $L_q$-centroid bodies $Z_q(\mu)$, $q\gr 1$, are defined through their support function
\begin{equation}\label{eq:intro-11}
h_{Z_q(\mu)}(y):= \|\langle \cdot ,y\rangle \|_{L_q(\mu)} = \left(\int_{{\mathbb R}^n}|\langle x,y\rangle|^qd\mu(x)\right)^{1/q},
\end{equation}
and have played a key role in the study of the distribution of linear functionals with respect to the
measure $\mu$. For every $1\ls q\ls n$ we obtain sharp upper bounds for the radius of random sections of $Z_q(\mu )$
of dimension proportional to $n$, thus extending a similar result of
Brazitikos and Stavrakakis which was established only for $q\in [1,\sqrt{n}]$.

\begin{theorem}\label{th:intro-Zq}Let $\mu $ be an isotropic log-concave probability measure on ${\mathbb R}^n$ and let $1\ls q\ls n$.
Then:
\begin{enumerate}
\item[{\rm (i)}] If $k=\gamma n$ for some $\gamma \in (0,1)$, then, with probability greater than $1-e^{-\overline{c}_4k}$, a random $F\in G_{n,n-k}$ satisfies
\begin{equation}\label{eq:intro-12}R(Z_q(\mu )\cap F)\ls \overline{c}_5(\gamma )\sqrt{q},\end{equation}
where $\overline{c}_4$ is an absolute constant and $\overline{c}_5(\gamma )=O(\gamma^{-2}\log^{5/2}(c/\gamma ))$ is a positive constant depending only on $\gamma $.
\item[{\rm (ii)}] With probability greater than $1-e^{-n}$, a random $U\in O(n)$ satisfies
\begin{equation}\label{eq:intro-13}Z_q(\mu )\cap U(Z_q(\mu ))\subseteq (\overline{c}_6\sqrt{q})\,B_2^n,\end{equation}
where $\overline{c}_6>0$ is an absolute constant.
\end{enumerate}
\end{theorem}

The method of proof is based on estimates (from \cite{EMilman-2014} and \cite{Giannopoulos-EMilman-2014}) for the Gelfand numbers of symmetric convex bodies in terms of their volumetric parameters; combining these general estimates with fundamental (known) properties of the family of the
centroid bodies $Z_q(\mu )$ of an isotropic log-concave probability measure $\mu $ we provide estimates for the {\sl minimal} radius of a $k$-codimensional
section of $Z_q(\mu )$. Then, we pass to bounds for the radius of random $k$-codimensional sections of $Z_q(\mu )$ using known results
from \cite{Giannopoulos-Milman-Tsolomitis-2005}, \cite{Vershynin-2006} and \cite{Litvak-Pajor-Tomczak-2006}. We conclude Section 4 with
a discussion of the same questions for the polar bodies $Z_q^{\circ }(\mu )$ of the centroid bodies $Z_q(\mu )$.

\smallskip

Using the same approach we study the diameter of random sections of convex bodies which have {\it maximal isotropic constant}. Set
\begin{equation}\label{eq:intro-22}L_n^{\prime }:= \max\{ L_K:K\ \hbox{is an isotropic symmetric convex body in}\ {\mathbb R}^n\}.\end{equation}
It is known that $L_n\ls cL_n^{\prime }$ for some absolute constant $c>0$ (see \cite[Chapter 3]{BGVV-book}). We prove the following:

\begin{theorem}\label{th:intro-max}Assume that $K$ is an isotropic symmetric convex body in ${\mathbb R}^n$ with $L_K=L_n^{\prime }$.
Then:
\begin{enumerate}
\item[{\rm (i)}] A random $F\in G_{n,n/2}$ satisfies
\begin{equation}\label{eq:max-1}R(K\cap F)\ls \overline{c}_7\sqrt{n}\end{equation}
and
\begin{equation}\label{eq:max-2}L_{K\cap F}\ls \overline{c}_8\end{equation} with probability greater than $1-e^{-\overline{c}_9n}$, where $\overline{c}_i>0$ are absolute constants.
\item[{\rm (ii)}] A random $U\in O(n)$ satisfies
\begin{equation}K\cap U(K)\subseteq (\overline{c}_{10}\sqrt{n})\,B_2^n,\end{equation} with probability greater than $1-e^{-n}$, where $\overline{c}_{10}>0$ is an absolute constant.
\end{enumerate}
\end{theorem}

The same arguments work if we assume that $K$ has almost maximal isotropic constant, i.e. $L_K\gr\beta L_n^{\prime }$ for some (absolute)
constant $\beta\in (0,1)$. We can obtain similar results, with the constants $\overline{c}_i$ now depending only on $\beta $. It should be noted
that Alonso-Guti\'{e}rrez, Bastero, Bernu\'{e}s and Paouris \cite{Alonso-Bastero-Bernues-Paouris-2010} have proved that every convex body $K$
has a section $K\cap F$ of dimension $n-k$ with isotropic constant
\begin{equation}L_{K\cap F}\ls c\sqrt{\frac{n}{k}}\log \Big (\frac{en}{k}\Big ).\end{equation}
For the proof of this result they considered an $\alpha $-regular $M$-position of $K$.
In Theorem \ref{th:intro-max} we consider convex bodies in the isotropic position and the estimates \eqref{eq:max-1} and \eqref{eq:max-2} hold
for a random subspace $F$.

\medskip

\noindent {\bf Radius of random projections of $L_q$-centroid bodies and sub-Gaussian subspaces of isotropic convex bodies.}
Let $K$ be a centered convex body of volume $1$ in $\mathbb R^n$. We say that a direction $\theta\in S^{n-1}$ is a
$\psi_{\alpha }$-direction (where $1\ls\alpha\ls 2$) for $K$ with constant $b>0$ if
\begin{equation}\|\langle\cdot ,\theta\rangle\|_{L_{\psi_{\alpha }}(K)}\ls b\|\langle\cdot, \theta\rangle\|_2,\end{equation}
where \begin{equation} \|\langle \cdot ,\theta\rangle \|_{L_{\psi_{\alpha}}(K)}:=\inf \left \{ t>0 :
\int_K\exp \big((|\langle x,\theta\rangle |/t)^\alpha\big) \,
dx\ls 2 \right\}.
\end{equation}Markov's inequality implies that if $K$ satisfies a
$\psi_\alpha$-estimate with constant $b$ in the direction of
$\theta$ then for all $t\gr 1$ we have $|\{x\in K : |\langle x,\theta
\rangle|\gr t \|\langle\cdot ,\theta\rangle\|_2\}|\ls
2e^{-t^a/b^\alpha}$. Conversely, one can check that tail estimates of this form
imply that $\theta $ is a $\psi_{\alpha }$-direction for $K$.

It is well-known that every $\theta \in S^{n-1}$ is a $\psi_1$-direction for $K$ with an absolute constant $C$. An open question is
if there exists an absolute constant $C>0$ such that every $K$ has at least one sub-Gaussian direction ($\psi_2$-direction) with constant $C$.
It was first proved by Klartag in \cite{Klartag-2007} that
for every centered convex body $K$ of volume $1$ in ${\mathbb R}^n$ there exists $\theta \in S^{n-1}$ such that
\begin{equation} |\{ x\in K: |\langle x, \theta \rangle | \gr
c t \|\langle\cdot, \theta\rangle\|_2 \}| \ls e^{-\frac{t^{2}}{[\log(t+1)]^{2a}}} \end{equation} for all
$t\gr 1$, where $a=3$ (equivalently, $\|\langle \cdot ,\theta\rangle\|_{L_{\psi_2}(K)}\ls C(\log n)^a\|\langle\cdot ,\theta\rangle\|_2$).
This estimate was later improved by Giannopoulos, Paouris and Valettas in \cite{Giannopoulos-Paouris-Valettas-2011}
and \cite{Giannopoulos-Paouris-Valettas-2012b} (see also \cite{Giannopoulos-Pajor-Paouris-2007})
 who showed that the body $\Psi_2(K)$ with support function $y\mapsto \|\langle \cdot ,y\rangle\|_{L_{\psi_2}(K)}$
has volume
\begin{equation}\label{eq:vol-ratio-psi-body-K}c_1\ls \left(\frac{|\Psi_{2}(K)|}{|Z_2(K)|}\right)^{1/n}\ls
c_2\sqrt{\log n}.\end{equation}
From \eqref{eq:vol-ratio-psi-body-K} it follows that there exists at least one
sub-Gaussian direction for $K$ with constant $b\ls C\sqrt{\log n}$.

Brazitikos and Hioni in \cite{Brazitikos-Hioni-2015} proved that if $K$ is isotropic then logarithmic bounds for $\|\langle \cdot ,\theta\rangle\|_{L_{\psi_2}(K)}$
hold true with probability polynomially close to $1$: For any $a>1$ one has
$$\|\langle \cdot ,\theta\rangle\|_{L_{\psi_2}(K)}\ls C(\log n)^{3/2}\max\left\{\sqrt{\log n},\sqrt{a}\right\}L_K$$
for all $\theta $ in a subset $\Theta_a $ of $S^{n-1}$ with $\sigma (\Theta_a )\gr 1-n^{-a}$, where $C>0$
is an absolute constant.

Here, we consider the question if one can have an estimate of this type {\sl for all} directions $\theta $ of a subspace $F\in G_{n,k}$
of dimension $k$ increasing to infinity with $n$. We say that $F\in G_{n,k}$ is a {\sl sub-Gaussian subspace} for $K$ with constant $b>0$
if
\begin{equation}\|\langle\cdot ,\theta\rangle\|_{L_{\psi_{\alpha }}(K)}\ls b\|\langle\cdot, \theta\rangle\|_2\end{equation}
for all $\theta \in S_F:=S^{n-1}\cap F$. In Section 5 we show that if $K$ is isotropic then a random subspace of dimension $(\log n)^4$
is sub-Gaussian with constant $b\simeq (\log n)^2$. More precisely, we prove the following.

\begin{theorem}\label{th:1.3}Let $K$ be an isotropic convex body in ${\mathbb R}^n$. If $k\simeq (\log n)^4$ then there exists a
subset $\Gamma $ of $G_{n,k}$ with $\nu_{n,k}(\Gamma )\gr 1-n^{-(\log n)^3}$ such that
\begin{equation}\|\langle \cdot ,\theta\rangle\|_{L_{\psi_2}(K)}\ls C(\log n)^2L_K\end{equation}
for all $F\in \Gamma $ and all $\theta\in S_F$, where $C>0$ is an absolute constant.
\end{theorem}

An essential ingredient of the proof is the good estimates on the radius of random projections of the $L_q$-centroid bodies $Z_q(K)$ of $K$,
 which follow from E.~Milman's sharp bounds on their mean width $w(Z_q(K))$ (see Theorem \ref{th:Emanuel2}).

\section{Notation and preliminaries}

We work in ${\mathbb R}^n$, which is equipped with a Euclidean structure $\langle\cdot ,\cdot\rangle $. We denote the corresponding
Euclidean norm by $\|\cdot \|_2$, and write $B_2^n$ for the Euclidean unit ball, and $S^{n-1}$ for the unit sphere. Volume is
denoted by $|\cdot |$. We write $\omega_n$ for the volume of $B_2^n$ and $\sigma $ for the rotationally invariant probability measure on
$S^{n-1}$. We also denote the Haar measure on $O(n)$ by $\nu $. The Grassmann manifold $G_{n,k}$ of $k$-dimensional subspaces of
${\mathbb R}^n$ is equipped with the Haar probability measure $\nu_{n,k}$. Let $k\ls n$ and $F\in G_{n,k}$. We will denote the
orthogonal projection from $\mathbb R^{n}$ onto $F$ by $P_F$. We also define $B_F=B_2^n\cap F$ and $S_F=S^{n-1}\cap
F$.

The letters $c,c^{\prime }, c_1, c_2$ etc. denote absolute positive constants whose value may change from line to line. Whenever we
write $a\simeq b$, we mean that there exist absolute constants $c_1,c_2>0$ such that $c_1a\ls b\ls c_2a$.  Also if $A,D\subseteq
\mathbb R^n$ we will write $A\simeq D$ if there exist absolute constants $c_1, c_2>0$ such that $c_{1}A\subseteq D \subseteq
c_{2}A$.

\medskip

\noindent \textbf{Convex bodies.} A convex body in ${\mathbb R}^n$ is a compact convex subset $A$ of
${\mathbb R}^n$ with nonempty interior. We say that $A$ is symmetric if $A=-A$. We say that $A$ is centered if
the center of mass of $A$ is at the origin, i.e.~$\int_A\langle
x,\theta\rangle \,d x=0$ for every $\theta\in S^{n-1}$.

The volume radius of $A$ is the quantity ${\rm vrad}(A)=\left (|A|/|B_2^n|\right )^{1/n}$.
Integration in polar coordinates shows that if the origin is an interior point of $A$ then the volume radius of $A$ can be expressed as
\begin{equation}\label{eq:not-1}{\rm vrad}(A)=\left (\int_{S^{n-1}}\|\theta \|_A^{-n}\,d\sigma (\theta )\right)^{1/n},\end{equation}
where $\|\theta \|_A=\min\{ t>0:\theta \in tA\}$. The radial function of $A$ is defined by $\rho_A(\theta )=\max\{ t>0:t\theta\in A\}$,
$\theta\in S^{n-1}$. The support
function of $A$ is defined by $h_A(y):=\max \bigl\{\langle x,y\rangle :x\in A\bigr\}$, and
the mean width of $A$ is the average
\begin{equation}\label{eq:not-2}w(A):=\int_{S^{n-1}}h_A(\theta )\,d\sigma (\theta )\end{equation}
of $h_A$ on $S^{n-1}$. The radius $R(A)$ of $A$ is the smallest $R>0$ such that $A\subseteq RB_2^n$.
For notational convenience we write $\overline{A}$ for
the homothetic image of volume $1$ of a convex body $A\subseteq
\mathbb R^n$, i.e. $\overline{A}:= |A|^{-1/n}A$.

The polar body $A^{\circ }$ of a convex body $A$ in ${\mathbb R}^n$ with $0\in {\rm int}(A)$ is defined by
\begin{equation}\label{eq:not-3}
A^{\circ}:=\bigl\{y\in {\mathbb R}^n: \langle x,y\rangle \ls 1\;\hbox{for all}\; x\in A\bigr\}.
\end{equation}
The Blaschke-Santal\'{o} inequality states that if $A$ is centered then $|A||A^{\circ }|\ls |B_2^n|^2$,
with equality if and only if $A$ is an ellipsoid.
The reverse Santal\'{o} inequality of J.~Bourgain and V.~Milman \cite{Bourgain-VMilman-1987} states that there exists an absolute constant $c>0$ such
that
\begin{equation}\label{eq:not-4}\left (|A||A^{\circ }|\right )^{1/n}\gr c/n\end{equation}
whenever $0\in {\rm int}(A)$.

For every centered convex body $A$ of volume $1$ in ${\mathbb R}^n$ and for every $q\in (-n,\infty )\setminus\{0\}$ we define
\begin{equation}I_q(A)=\left (\int_A\| x\|_2^qdx\right )^{1/q}.\end{equation}
As a consequence of Borell's lemma (see \cite[Chapter 1]{BGVV-book}) one has
\begin{equation}I_q(A)\ls c_1 q I_2(A)\end{equation} for all $q\gr 2$.

For basic facts from the Brunn-Minkowski theory and the asymptotic theory of convex bodies we refer to the books \cite{Schneider-book} and \cite{AGA-book} respectively.

\smallskip

\noindent \textbf{Log-concave probability measures.}
Let $\mu $ be a log-concave probability measure on ${\mathbb R}^n$. The density of $\mu $ is denoted by $f_{\mu}$. We say that $\mu $
is centered and we write $\textrm{bar}(\mu )=0$ if, for all $\theta\in S^{n-1}$,
\begin{equation}\label{eq:not-5}
\int_{\mathbb R^n} \langle x, \theta \rangle d\mu(x) = \int_{\mathbb
R^n} \langle x, \theta \rangle f_{\mu}(x) dx = 0.
\end{equation}
The isotropic constant of $\mu $ is defined by
\begin{equation}\label{eq:definition-isotropic}
L_{\mu }:=\left (\frac{\sup_{x\in {\mathbb R}^n} f_{\mu} (x)}{\int_{{\mathbb
R}^n}f_{\mu}(x)dx}\right )^{\frac{1}{n}} [\det \textrm{Cov}(\mu)]^{\frac{1}{2n}},\end{equation} where
$\textrm{Cov}(\mu)$ is the covariance matrix of $\mu$ with entries
\begin{equation}\label{eq:not-6}\textrm{Cov}(\mu )_{ij}:=\frac{\int_{{\mathbb R}^n}x_ix_j f_{\mu}
(x)\,dx}{\int_{{\mathbb R}^n} f_{\mu} (x)\,dx}-\frac{\int_{{\mathbb
R}^n}x_i f_{\mu} (x)\,dx}{\int_{{\mathbb R}^n} f_{\mu}
(x)\,dx}\frac{\int_{{\mathbb R}^n}x_j f_{\mu}
(x)\,dx}{\int_{{\mathbb R}^n} f_{\mu} (x)\,dx}.\end{equation} We say
that a log-concave probability measure $\mu $ on ${\mathbb R}^n$
is isotropic if $\textrm{bar}(\mu )=0$ and $\textrm{Cov}(\mu )$ is the identity matrix.
Note that a centered convex body $K$ of volume $1$ in ${\mathbb R}^n$ is isotropic,
i.e.~it satisfies (\ref{eq:intro-1}),
if and only if the log-concave probability measure $\mu_K$ with density
$x\mapsto L_K^n\mathbf{1}_{K/L_K}(x)$ is isotropic. Note that for every log-concave measure $\mu $
on ${\mathbb R}^n$ one has
\begin{equation}\label{eq:Lmu}L_{\mu }\ls \kappa L_n,\end{equation}
where $\kappa >0$ is an absolute constant (a proof can be found in \cite[Proposition 2.5.12]{BGVV-book}).

We will use the following sharp result on the growth of $I_q(K)$, where $K$ is an isotropic
convex body in ${\mathbb R}^n$, proved by Paouris in \cite{Paouris-GAFA} and \cite{Paouris-TAMS}.

\begin{theorem}[Paouris]\label{th:grigoris}
There exists an absolute constant $\delta >0$ with the following property: if $K$ is an isotropic convex body in $\mathbb
R^n$, then
\begin{equation}\label{eq:constant-Iq}\frac{1}{\delta }\sqrt{n}L_K=\frac{1}{\delta }I_2(K)\ls I_{-q}(K)\ls I_q(K)\ls \delta I_2(K)=\delta\sqrt{n}L_K\end{equation} for every
$1\ls q\ls \sqrt{n}$.
\end{theorem}

For every $q\gr 1$ and every $y \in {\mathbb R}^n$ we set
\begin{equation}\label{Zq-def}h_{Z_q(\mu )}(y)= \left(\int_{{\mathbb R}^n} |\langle x,y\rangle|^{q}d\mu (x) \right)^{1/q}.\end{equation}
The $L_q$-centroid body $Z_q(\mu )$ of $\mu $ is the symmetric convex body with support function
$h_{Z_{q}(\mu )}$. Note that $\mu $ is isotropic if and only if it is centered and $Z_{2}(\mu )=
B_2^n$. If $K$ is an isotropic convex body in ${\mathbb R}^n$ we define $Z_q(K)=L_KZ_q(\mu_K)$.
From H\"{o}lder's inequality it follows that $Z_1(K)\subseteq Z_p(K)\subseteq Z_q(K)\subseteq Z_{\infty }(K)$ for
all $1\ls p\ls q\ls \infty $, where $Z_{\infty }(K)={\rm conv}\{K,-K\}$.
Using Borell's lemma, one can check that
\begin{equation}\label{eq:Zq-inclusions} Z_q(K)\subseteq c_1\frac{q}{p}Z_p(K)\end{equation}
for all $1\ls p<q$. In particular, if $K$ is isotropic, then
$R(Z_q(K))\ls c_1qL_K$. One can also check that if $K$ is
centered, then $Z_q(K)\supseteq c_2Z_{\infty }(K)$ for all $q\gr n$.

It was shown by Paouris \cite{Paouris-GAFA} that if $1\ls q\ls\sqrt{n}$ then
\begin{equation}\label{eq:wZq-small}
w\bigl(Z_q(\mu)\bigr)\simeq \sqrt{q},
\end{equation}
and that for all $1\ls q\ls n$,
\begin{equation}
{\rm vrad}(Z_q(\mu)) \ls c_1\sqrt{q}.
\end{equation}
Conversely, it was shown by B.~Klartag and E.~Milman in \cite{Klartag-EMilman-2012} that if $1\ls q\ls\sqrt{n}$ then
\begin{equation}\label{eq:low-volume-Zq}
{\rm vrad}(Z_q(\mu))\gr c_2\sqrt{q}.
\end{equation}
This determines the volume radius of $Z_q(\mu )$ for all $1 \ls q\ls\sqrt{n}$. For larger values of $q$ one can still use the lower bound:
\begin{equation}\label{eq:3}
{\rm vrad}(Z_q(\mu)) \gr c_2\sqrt{q}\, L_{\mu}^{-1} ,
\end{equation}
obtained by Lutwak, Yang and Zhang in \cite{Lutwak-Yang-Zhang-2000} for convex bodies and extended by Paouris and Pivovarov
in \cite{Paouris-Pivovarov-2012} to the class of log-concave probability measures.

Let $\mu $ be a probability measure on ${\mathbb R}^n$ with density $f_{\mu }$ with respect
to the Lebesgue measure. For every $1\ls k\ls n-1$ and every
$E\in G_{n,k}$, the marginal of $\mu$ with respect to $E$ is the probability
measure $\pi_E(\mu )$ on $E$, with density
\begin{equation}\label{definitionmarginal}f_{\pi_E(\mu )}(x)= \int_{x+
E^{\perp}} f_{\mu }(y) dy.
\end{equation}
It is easily checked that if $\mu $ is centered, isotropic or log-concave, then $\pi_E(\mu )$ is also centered, isotropic or
log-concave, respectively. A very useful observation is that:
\begin{equation}
P_F\bigl(Z_q(\mu )\bigr) = Z_q\bigl(\pi_F(\mu )\bigr)
\end{equation}
for every $1\ls k\ls n-1$ and every $F\in G_{n,n-k}$.

If $\mu$ is a centered log-concave probability measure on $\mathbb R^n$ then for every $p>0$ we define
\begin{equation}\label{eq:not-7}
K_p(\mu):=K_p(f_\mu)=\left\{x : \int_0^\infty r^{p-1}f_\mu(rx)\,
dr\gr \frac{f_{\mu }(0)}{p} \right\}.
\end{equation}
From the definition it follows that $K_p(\mu )$ is a star body with radial function
\begin{equation}\label{eq:not-8}
\rho_{K_p(\mu )}(x)=\left (\frac{1}{f_{\mu }(0)}\int_0^{\infty
}pr^{p-1}f_{\mu }(rx)\,dr\right )^{1/p}\end{equation} for $x\neq 0$.
The bodies $K_p(\mu )$ were introduced in \cite{Ball-1988} by K. Ball who showed that if $\mu $ is log-concave
then, for every $p>0$, $K_p(\mu )$ is a convex body.

If $K$ is isotropic then for every $1\ls k\ls n-1$ and $F\in G_{n,n-k}$,
the body $\overline{K_{k+1}}(\pi_{F^{\perp }}(\mu_{K}))$ satisfies
\begin{equation}\label{eq:not-9}
|K\cap F|^{1/k} \simeq
\frac{L_{\overline{K_{k+1}}(\pi_{F^{\perp }}(\mu_{K}))}}{L_K}.
\end{equation}
For more information on isotropic convex bodies and log-concave measures see \cite{BGVV-book}.

\section{Random sections of isotropic convex bodies}

The proof of Theorem \ref{th:intro-2} is based on Lemma \ref{lem:sections-1} and Lemma \ref{lem:sections-2} below. They exploit
some ideas of Klartag from \cite{Klartag-2004}.

\begin{lemma}\label{lem:sections-1}Let $K$ be an isotropic convex body in ${\mathbb R}^n$. For every $1\ls k\ls n-1$
there exists a subset ${\cal A}:={\cal A}(n,k)$ of $G_{n,n-k}$ with $\nu_{n,n-k}({\cal A})\gr 1-e^{-\sqrt{n}}$ that has the following property:
for every $F\in {\cal A}$,
\begin{equation}\label{eq:lem-sections-1}|\{x\in K\cap F:\|x\|_2
\gr c_1\sqrt{n}L_K\}|\ls e^{-(k+\sqrt{n})}|K\cap F|,
\end{equation}
where $c_1>0$ is an absolute constant.
\end{lemma}

\noindent {\it Proof.} Integration in polar coordinates shows that for all $q>0$ 
\begin{equation}\label{eq:sections-1}\int_{G_{n,n-k}}\int_{K\cap F}\|x\|_2^{k+q}dx\,d\nu_{n,n-k}(F)
=\frac{(n-k)\omega_{n-k}}{n\omega_n}\int_K\|x\|_2^qdx =\frac{(n-k)\omega_{n-k}}{n\omega_n}I_q^q(K),\end{equation}
and an application of Markov's inequality shows that a random $F\in G_{n,n-k}$ satisfies
\begin{equation}\label{eq:sections-2}\int_{K\cap F}\|x\|_2^{k+q}dx\ls \frac{(n-k)\omega_{n-k}}{n\omega_n}(eI_q(K))^q\end{equation}
with probability greater than $1-e^{-q}$.

Fix a subspace $F\in G_{n,n-k}$ which satisfies \eqref{eq:sections-2}. From \eqref{eq:not-9} we have
\begin{equation}\label{eq:sections-3}|K\cap F|^{1/k} \gr c_2
\frac{L_{\overline{K_{k+1}}(\pi_{F^{\perp }}(\mu_{K}))}}{L_K}\gr \frac{c_3}{L_K}
\end{equation}
where $c_2,c_3>0$ are absolute constants. A simple computation shows that
\begin{equation}\frac{(n-k)\omega_{n-k}}{n\omega_n}\ls (c_4\sqrt{n})^k\end{equation}
for an absolute constant $c_4>0$. Using also \eqref{eq:constant-Iq} with $q=\sqrt{n}$ we get
\begin{align}\label{eq:sections-4}\frac{1}{|K\cap F|}\int_{K\cap F}\|x\|_2^{k+\sqrt{n}}dx &\ls  \frac{1}{|K\cap F|}\,\frac{(n-k)\omega_{n-k}}{n\omega_n}(eI_{\sqrt{n}}(K))^{\sqrt{n}}\\
\nonumber &\ls (c_5L_K)^k(c_4\sqrt{n})^k(e\delta\sqrt{n}L_K)^{\sqrt{n}}\ls (c_6\sqrt{n}L_K)^{k+\sqrt{n}},
\end{align}
where $c_6>0$ is an absolute constant. It follows that
\begin{equation}\label{eq:sections-5}|\{x\in K\cap F:\|x\|_2\gr ec_6\sqrt{n}L_K\}|\ls e^{-(k+\sqrt{n})}|K\cap F|.
\end{equation}
and the lemma is proved with $c_1=ec_6$. \prend

\medskip

The next lemma comes from \cite{Klartag-2004}.

\begin{lemma}[Klartag]\label{lem:sections-2}Let $A$ be a symmetric convex body in
$\mathbb R^m$. Then, for any $0<\varepsilon <1$ we have
\begin{equation}\label{eq:lem-sections-2}|\{x\in A
: \|x\|_2\gr \varepsilon R(A)\}|\gr \frac{1}{2}(1-\varepsilon)^m|A|.\end{equation}
\end{lemma}

\noindent {\it Proof.} Let $x_0\in A$ such that $\|x_0\|_2=R(A)$ and define $v=x_0/\|x_0\|_2$. We consider the set $A^+$
defined as
\begin{equation}A^+:=\{x\in A : \langle x,v\rangle \gr 0\}.\end{equation}
Since $A$ is symmetric, we have $|A^+|= |A|/2$. Note that
\begin{equation}\{x\in A :\|x\|_2\gr \varepsilon R(A)\}\supseteq \varepsilon x_0+
(1-\varepsilon)A^+.\end{equation} Therefore,
\begin{equation}|\{x\in A: \|x\|_2\gr \varepsilon
R(A)\}|\gr |\varepsilon x_0+(1-\varepsilon) A^+|=(1-\varepsilon)^m|A^+|= \frac{1}{2}(1-\varepsilon)^m|A|,\end{equation}
as claimed. \prend

\medskip

\noindent {\bf Proof of Theorem \ref{th:intro-2}.} Let $K$ be an isotropic symmetric convex body in ${\mathbb R}^n$. Applying
Lemma \ref{lem:sections-1} we find a subset ${\cal A}$ of $G_{n,n-k}$ with $\nu_{n,n-k}({\cal A})\gr 1-e^{-\sqrt{n}}$
such that, for every $F\in {\cal A}$,
\begin{equation}\label{eq:final-1}|\{x\in K\cap F:\|x\|_2
\gr c_1\sqrt{n}L_K\}|\ls e^{-(k+\sqrt{n})}|K\cap F|.
\end{equation}
We distinguish two cases:

\smallskip

\noindent {\it Case 1.} If $k>n/3$ then choosing $\varepsilon_0 =1-e^{-\frac{1}{3}}$ we get
\begin{equation}\frac{1}{2}(1-\varepsilon_0)^{n-k}|K\cap F|=\frac{1}{2}e^{-\frac{n-k}{3}}|K\cap F|>e^{-\frac{n-k}{3}-1}|K\cap F|
>e^{-(k+\sqrt{n})}|K\cap F|,\end{equation}
because $k+\sqrt{n}>\frac{n-k}{3}+1$. By Lemma \ref{lem:sections-2} and (\ref{eq:final-1})
we get that
\begin{equation}|\{x\in K\cap F: \|x\|_2\gr \varepsilon_0 R(K\cap F)\}|>|\{x\in K\cap F:\|x\|_2
\gr c_1\sqrt{n}L_K\}|,\end{equation}
therefore
\begin{equation}R(K\cap F)<c_2\sqrt{n}L_K,\end{equation}
where $c_2=\varepsilon_0^{-1}c_1>0$ is an absolute constant.

\smallskip

\noindent {\it Case 2.} If $k\ls n/3$ then we choose $\varepsilon_1 =\frac{k+\sqrt{n}}{6(n-k)}$. Note that $\varepsilon_1<1/2$.
Using the inequality $1-t>e^{-2t}$ on $(0,1/2)$ we get
\begin{equation}\frac{1}{2}(1-\varepsilon_1)^{n-k}|K\cap F|=\frac{1}{2}\left (1-\frac{k+\sqrt{n}}{6(n-k)}\right )^{n-k}|K\cap F|
>e^{-\frac{k+\sqrt{n}}{3}-1}|K\cap F|>e^{-(k+\sqrt{n})}|K\cap F|,\end{equation}
because $\frac{2(k+\sqrt{n})}{3}>1$. By Lemma \ref{lem:sections-2} this implies that
\begin{equation}|\{x\in K\cap F: \|x\|_2\gr \varepsilon_1 R(K\cap F)\}|>|\{x\in K\cap F:\|x\|_2
\gr c_1\sqrt{n}L_K\}|,\end{equation}
therefore
\begin{equation}\varepsilon_1 R(K\cap F)< c_1\sqrt{n}L_K,\end{equation}
which, by the choice of $\varepsilon_1$ becomes
\begin{equation}R(K\cap F)< \frac{c_3n}{\max\{ k,\sqrt{n}\}}\,\sqrt{n}L_K\end{equation}
for some absolute constant $c_3>0$. This completes the proof of the theorem (with a probability estimate $1-e^{-\sqrt{n}}$ for all $1\ls k\ls n-1$).
\prend

\begin{remark}\label{rem:3-2}\rm It is possible to improve the probability estimate $1-e^{-\sqrt{n}}$ in the range $k\gr \gamma n$, for any $\gamma\in (1/\sqrt{n},1)$.
This can be done with the help of known results that demonstrate the fact
that the existence of one $s$-dimensional section with radius $r$ implies that random $m$-dimensional sections, where $m<s$, have
radius of ``the same order". This was first observed in \cite{Giannopoulos-Milman-Tsolomitis-2005}, \cite{Vershynin-2006} and, soon after,
in \cite{Litvak-Pajor-Tomczak-2006}. Let us recall this last statement.
\begin{quote}{\sl Let $A$ be a symmetric convex body in ${\mathbb R}^n$ and let
$1\ls s<m\ls n-1$. If $R(A\cap F)\ls r $ for some $F\in G_{n,m}$ then a random subspace $E\in G_{n,s}$ satisfies
\begin{equation}\label{eq:LPT}R(A\cap E)\ls r\,\Big ( \frac{c_2n}{n-m}\Big )^{\frac{n-s}{2(m-s)}}\end{equation}
with probability greater than $1-2e^{-(n-s)/2}$, where $c_2>0$ is an absolute constant.}
\end{quote}
We apply this result as follows. Let $k=\gamma n\gr\sqrt{n}$ and set $t=\delta n$, where $\delta \simeq\gamma /\log (1+1/\gamma )$. From the proof of Theorem \ref{th:intro-2} we know that there exists $E\in G_{n,n-t}$ such that
\begin{equation}R(K\cap E)\ls \frac{c_1n}{t}\sqrt{n}L_K,\end{equation}
where $c_1>0$ is an absolute constant. Applying \eqref{eq:LPT} with $s=n-k$ and $m=n-t$ we see that
a random subspace $F\in G_{n,n-k}$ satisfies 
\begin{equation}R(K\cap F)\ls \left (\frac{c_2}{\delta }\right )^{\frac{3}{2}}\,R(K\cap E)= c_3(\gamma )\sqrt{n}L_K\end{equation}
with probability greater than $1-2e^{-k/2}$, where $c_3(\gamma )=O((\gamma^{-1}\log (1+1/\gamma ))^{\frac{3}{2}})$.
\end{remark}

\begin{remark}\label{rem:3-3}\rm It is also possible to give lower bounds of the order of $\sqrt{n}L_K$
for the diameter of $(n-k)$-dimensional sections, provided that the codimension $k$ is small. Integration in polar coordinates shows that
\begin{equation}\int_K\|x\|_2^{-q}dx=\frac{n\omega_n}{(n-k)\omega_{n-k}}\int_{G_{n,n-k}}\int_{K\cap F}\|x\|_2^{k-q}dx\,d\nu_{n,n-k}(F)\end{equation}
for every $1\ls k\ls n-1$ and every $0<q<n$. It follows that
\begin{equation}\label{eq:sections-lower-1}\int_{G_{n,n-k}}\int_{K\cap F}\|x\|_2^{k-q}dx\,d\nu_{n,n-k}(F)=\frac{(n-k)\omega_{n-k}}{n\omega_n}I_{-q}^{-q}(K),\end{equation}
and an application of Markov's inequality shows that a random $F\in G_{n,n-k}$ satisfies
\begin{equation}\label{eq:sections-lower-2}\int_{K\cap F}\|x\|_2^{k-q}dx\ls \frac{(n-k)\omega_{n-k}}{n\omega_n}(e/I_{-q}(K))^q\end{equation}
with probability greater than $1-e^{-q}$. Assuming that $q>k$, for any $F\in G_{n,n-k}$ satisfying \eqref{eq:sections-lower-2} we have
\begin{equation}\label{eq:sections-lower-3}|K\cap F|\,R(K\cap F)^{k-q}\ls \int_{K\cap F}\|x\|_2^{k-q}dx\ls \frac{(n-k)\omega_{n-k}}{n\omega_n}(e/I_{-q}(K))^q,\end{equation}
which implies
\begin{equation}\label{eq:sections-lower-4}R(K\cap F)\gr \left (\frac{n\omega_n}{(n-k)\omega_{n-k}}\right )^{\frac{1}{q-k}}|K\cap F|^{\frac{1}{q-k}}\left (\frac{I_{-q}(K)}{e}\right )^{\frac{q}{q-k}}\gr \left (\frac{c_1}{\sqrt{n}L_K}\right )^{\frac{k}{q-k}}\left (c_2I_{-q}(K)\right )^{\frac{q}{q-k}}.\end{equation}
If $k\ls\sqrt{n}$ then we may choose $q=2\sqrt{n}$ and use the fact that $I_{-2\sqrt{n}}(K)\gr c_3\sqrt{n}L_K$ by Theorem \ref{th:grigoris}, to get:
\end{remark}

\begin{proposition}\label{prop:sections-3}Let $K$ be an isotropic convex body in ${\mathbb R}^n$. For every $1\ls k\ls\sqrt{n}$
there exists a subset ${\cal A}$ of $G_{n,n-k}$ with $\nu_{n,n-k}({\cal A})\gr 1-e^{-\sqrt{n}}$
such that, for every $F\in {\cal A}$,
\begin{equation}\label{eq:prop-lower-1}R(K\cap F)\gr c\sqrt{n}L_K,\end{equation}
where $c>0$ is an absolute constant.
\end{proposition}

\begin{remark}\label{rem:3-1}\rm Choosing $k=\lfloor n/2\rfloor $ in Theorem \ref{th:intro-2} we see that if $K$ is an isotropic symmetric convex body in ${\mathbb R}^n$
then a subspace $F\in G_{n,\lceil n/2\rceil }$ satisfies
\begin{equation}\label{eq:rotations-1}R(K\cap F)\ls c_1\sqrt{n}\,L_K\end{equation}
with probability greater than $1-2\exp (-c_2n)$, where $c_1, c_2>0$ are absolute
constants. A standard argument that goes back to Krivine (see \cite[Proposition 8.6.2]{AGA-book})
shows that there exists $U\in O(n)$ such that
\begin{equation}\label{eq:rotations-2}K\cap U(K)\subseteq (c_3\sqrt{n}L_K)\,B_2^n,\end{equation}
where $c_3>0$ is an absolute constant. In fact, one can prove an analogue of \eqref{eq:rotations-2} for a random $U\in O(n)$ using a result of Vershynin and Rudelson
(see \cite[Theorem 1.1]{Vershynin-2006}): There exist absolute constants $\gamma_0\in (0,1/2)$ and $c_1>0$ with the following property: if $A$ and $D$ are two symmetric convex bodies in ${\mathbb R}^n$ which have sections of dimensions at least $k$ and $n-2\gamma_0k$ respectively whose radius is bounded by $1$, then a random $U\in O(n)$
satisfies
\begin{equation}\label{eq:Rud-Ver}R(A\cap U(D))\ls c_1^{n/k}\end{equation} with probability greater than $1-e^{-n}$. As an application, setting $D=A$ and $k=n/2$ one has the following
(see \cite{Brazitikos-Stavrakakis-2014}). If
\begin{equation}r_A:=\min\{ R(A\cap F):{\rm dim}(F)=\lceil (1-\gamma_0)n\rceil \}\end{equation}
then $R(A\cap U(A))\ls c_2r_A$ with probability greater than $1-e^{-n}$ with respect to $U\in O(n)$.

Choosing $k=\lfloor \gamma_0n/2\rfloor $ in Theorem \ref{th:intro-2} we see that if $K$ is an isotropic symmetric convex body in ${\mathbb R}^n$
then
\begin{equation}r_K\ls c_4\sqrt{n}L_K\end{equation}
for some absolute constant $c_4>0$. This gives that a random $U\in O(n)$ satisfies
\begin{equation}K\cap U(K)\subseteq (c_5\sqrt{n}L_K)\,B_2^n,\end{equation} with probability greater than $1-e^{-n}$, where $c_5>0$ is an absolute constant.
\end{remark}

\section{Minimal and random sections of the centroid bodies of isotropic log-concave measures}

In this section we discuss the case of the $L_q$-centroid bodies $Z_q(\mu )$ of an isotropic log-concave probability
measure $\mu $ on ${\mathbb R}^n$. Our method will be different from the one in the previous section.

In view of \eqref{eq:LPT} we can give an upper bound for the radius of a
random $k$-codimensional section of a symmetric convex body $A$ in ${\mathbb R}^n$ if we are able to give an upper bound
for the radius of {\it some} $t$-codimensional section of $A$, where $t\ll k$. This leads us to the study of the
Gelfand numbers $c_t(A)$, which are defined by
\begin{equation}c_t(A)=\min\{R(A\cap F) : F \in G_{n,n-t}\}\end{equation}
for every $t = 0,\ldots,n-1$.
It was proved in \cite{Giannopoulos-EMilman-2014} that if $A$ is a symmetric convex body in ${\mathbb R}^n$
then, for any $t=1,\ldots,\lfloor n/2 \rfloor$ there exists $F \in G_{n,n-2t}$ such that
\begin{equation} \label{eq:G-EM-1}
A\cap F \subseteq c_1\frac{n}{t} \log\Big(e + \frac{n}{t}\Big) w_{t}(A) B_2^n \cap F,
\end{equation} where
\begin{equation}w_t(A) := \sup \{{\rm vrad}(A\cap E): E\in G_{n,t}\}.\end{equation}
In other words,
\begin{equation} \label{eq:G-EM-2}
c_{2t}(A)\ls c_1\frac{n}{t} \log\Big(e + \frac{n}{t}\Big) w_t(A).
\end{equation}
This is a refinement of a result of V.~Milman and G.~Pisier from \cite{VMilman-Pisier-1987}, where a similar estimate was
obtained, with the parameter $w_t(A)$ replaced by (the larger one)
\begin{equation}\label{eq:MP-parameter}v_t(A):= \sup \{{\rm vrad}(P_E(A)): E\in G_{n,t}\}.\end{equation}
We shall apply this method to the bodies $Z_q(\mu )$. The main additional
ingredient is the next fact, which combines results of Paouris and Klartag (see \cite{EMilman-2014} or \cite[Chapter 5]{BGVV-book}
for precise references):

\begin{theorem}\label{th:Emanuel}Let $\mu $ be a centered log-concave probability
measure on ${\mathbb R}^n$. Then, for all $1\ls t\ls n$ and $q\gr 1$ we have
\begin{equation}v_t(Z_q(\mu ))= \sup \{{\rm vrad}(P_E(Z_q(\mu ))): E\in G_{n,t}\}
\ls c_0\sqrt{\frac{q}{t}}\max\{ \sqrt{q},\sqrt{t}\}\max_{E\in G_{n,t}}\det\,{\rm Cov}(\pi_E(\mu ))^{\frac{1}{2t}},
\end{equation}
where $c_0>0$ is an absolute constant.
\end{theorem}

We apply Theorem \ref{th:Emanuel} as follows: for every $1\ls t\ls n/2$ and every $E\in G_{n,t}$ we have that $\pi_E(\mu )$ is isotropic, and
hence $\det\,{\rm Cov}(\pi_E(\mu ))^{\frac{1}{2t}}=1$. Then,
\begin{equation}w_t(Z_q(\mu ))\ls v_t(Z_q(\mu ))\ls c_0\sqrt{\frac{q}{t}}\max\{ \sqrt{q},\sqrt{t}\}.\end{equation}
From \eqref{eq:G-EM-2} we get

\begin{lemma}\label{lem:Zq-1}Let $\mu $ be an isotropic log-concave probability measure on ${\mathbb R}^n$
and let $1\ls t\ls \lfloor n/2 \rfloor$ and $1\ls q\ls n$. Then,
\begin{equation}\label{eq:ct-Zq}c_{2t}(Z_q(\mu ))\ls c_2\frac{n}{t}\log\left (e+\frac{n}{t}\right )\sqrt{\frac{q}{t}}\max\{ \sqrt{q},\sqrt{t}\},\end{equation}
where $c_2>0$ is an absolute constant. \prend
\end{lemma}

Let $k\gr 4$ and let $t<k/2$. From Lemma \ref{lem:Zq-1} we know that there exists $E\in G_{n,n-2t}$ such that
\begin{equation}R(Z_q(\mu )\cap E)\ls c_2\frac{n}{t}\log\left (e+\frac{n}{t}\right )\sqrt{\frac{q}{t}}\max\{ \sqrt{q},\sqrt{t}\},\end{equation}
where $c_2>0$ is an absolute constant. Applying \eqref{eq:LPT} with $s=n-k$ and $m=n-2t$ we see that
a random subspace $F\in G_{n,n-k}$ satisfies
\begin{equation}R(Z_q(\mu )\cap F)\ls \left (\frac{c_2n}{t}\right )^{\frac{k}{2(k-2t)}}\,R(Z_q(\mu )\cap E)\ls \left (\frac{c_3n}{t}\right )^{\frac{3}{2}+\frac{t}{k-2t}}\log\left (e+\frac{n}{t}\right )
\sqrt{\frac{q}{t}}\max\{ \sqrt{q},\sqrt{k}\}\end{equation}
with probability greater than $1-2e^{-k/2}$, where $c_3>0$ is an absolute constant. In particular, if $k=\gamma n$ we can choose
$t=\gamma n /\log (c/\gamma )$, for $c>e^2$, to get the following.

\begin{theorem}\label{th:Zq-sections}Let $\mu $ be an isotropic log-concave probability measure on ${\mathbb R}^n$
and let $\gamma \in (0,1)$ and $1\ls q\ls n$. If $k\gr \gamma n$ then a random subspace $F\in G_{n,n-k}$ satisfies
\begin{equation}R(Z_q(\mu )\cap F)\ls c(\gamma )\sqrt{q}\end{equation}
with probability greater than $1-2e^{-\gamma n/2}$, where $c(\gamma )=O(\gamma^{-2}\log^{5/2}(c/\gamma ))$ is a positive constant depending only on $\gamma $.
\end{theorem}

Next, we apply \eqref{eq:Rud-Ver}: choosing $t=\gamma_0n/2$ in \eqref{eq:ct-Zq} we see that
\begin{equation}r_{Z_q(\mu )}=c_{\gamma_0n}(Z_q(\mu ))\ls c_4\sqrt{q} \end{equation}
for every $1\ls q\ls n$, where $c_4=c_4(\gamma_0)>0$ is an absolute constant. Therefore, we have:

\begin{theorem}\label{th:Zq-rotations}Let $\mu $ be an isotropic log-concave probability measure on ${\mathbb R}^n$ and let $1\ls q\ls n$.
Then, a random $U\in O(n)$ satisfies
\begin{equation}Z_q(\mu )\cap U(Z_q(\mu ))\subseteq (c\sqrt{q})\,B_2^n,\end{equation}
with probability greater than $1-e^{-n}$, where $c>0$ is an absolute constant.
\end{theorem}

Note that Theorem \ref{th:intro-Zq} summarizes the contents of Theorem \ref{th:Zq-sections} and Theorem \ref{th:Zq-rotations}.

\begin{remark}\label{rem:Zq-polars}\rm We can study the same question for the polar body $Z_q^{\circ }(\mu )$ of $Z_q(\mu )$. Note that
\begin{equation}w_t(Z_q^{\circ }(\mu )) := \sup \{{\rm vrad}(Z_q^{\circ }(\mu )\cap E): E\in G_{n,t}\}
\simeq [\inf \{{\rm vrad}(P_E(Z_q(\mu ))): E\in G_{n,t}\}]^{-1}\end{equation}
by duality and by the Bourgain-Milman inequality. For any $1\ls t\ls n-1$ and any symmetric convex body $A$ in ${\mathbb R}^n$ define
\begin{equation}v_t^-(A)=\inf \{{\rm vrad}(P_E(A)): E\in G_{n,t}\}.\end{equation}
In the case $A=Z_q(\mu )$ this parameter has been studied in \cite{Giannopoulos-EMilman-2014}:

\begin{lemma}\label{lem:vk-Zq}
Let $\mu $ be an isotropic log-concave probability measure on ${\mathbb R}^n$. For any $q\gr 1$ and $1\ls k\ls n-1$ we have:
\begin{equation}\label{eq:vk-Zq-1}v_k^{-}(Z_q(\mu)) \gr c_1 \sqrt{\min(q, \sqrt{k})}.\end{equation}
If we assume that $\sup_nL_n\ls\alpha $ then we have
\begin{equation}\label{eq:vk-Zq-2}v_k^{-}(Z_q(\mu)) \gr \frac{c_2}{\alpha }\sqrt{\min(q, k)}\end{equation}
\end{lemma}

These estimates are leading to the next bounds on the minimal radius of a $k$-codimensional section of $Z_q^{\circ }(\mu )$. The following theorem is also from \cite{Giannopoulos-EMilman-2014}.

\begin{theorem}\label{thm:Rqk-variant}
Let $\mu $ be an isotropic log-concave probability measure on ${\mathbb R}^n$. For any $q\gr 1$ and $1\ls k\ls n-1$ we have:

\smallskip

{\rm (i)} There exists $F \in G_{n,n-k}$ such that:
\begin{equation}\label{eq:Rqk-variant-1}P_F(Z_q(\mu)) \supseteq \frac{1}{R_{k,q}}B_2^n\cap F\quad\textrm{and hence}\quad R(Z_q^{\circ }(\mu )\cap F)\ls R_{k,q},\end{equation}
where
\begin{equation}\label{eq:Rqk-variant-2}R_{k,q} = \min\left\{ 1, c_3\frac{1}{\min(q^{1/2} , k^{1/4})} \frac{n}{k} \log\left( e+ \frac{n}{k}\right) \right\} .
\end{equation}

{\rm (ii)} If we assume that $\sup_nL_n\ls\alpha $ then there exists $F \in G_{n,n-k}$ such that:
\begin{equation}\label{eq:Rqk-variant-3}P_F(Z_q(\mu)) \supseteq \frac{1}{R_{k,q,\alpha }}B_2^n\cap F\quad\textrm{and hence}\quad R(Z_q^{\circ }(\mu )\cap F)\ls R_{k,q,\alpha },\end{equation}
where
\begin{equation}\label{eq:Rqk-variant-4}R_{k,q,\alpha } = \min \left\{ 1, c_4\alpha\frac{1}{\sqrt{\min(q,k)}} \frac{n}{k} \log\left( e+ \frac{n}{k}\right) \right\}.
\end{equation}
\end{theorem}

Assuming that $q\ls\sqrt{n}$ and choosing $k=\gamma_0n$ we see from \eqref{eq:Rqk-variant-1} and \eqref{eq:Rqk-variant-2} that
\begin{equation}c_{\gamma_0n}(Z_q^{\circ }(\mu ))\ls c_1(\gamma_0)\frac{1}{\sqrt{q}}\end{equation}
where $c_1(\gamma_0)>0$ is an absolute constant. Then, we apply \eqref{eq:LPT} with $s=n/2$ and $m=(1-\gamma_0)n$
to get that a random subspace $E\in G_{n,n/2}$ satisfies
\begin{equation}R(Z_q^{\circ }(\mu )\cap E)\ls c_3\cdot c_{\gamma_0n}(Z_q^{\circ }(\mu ))\ls c_2(\gamma_0)\frac{1}{\sqrt{q}}\end{equation}
with probability greater than $1-2e^{-n/4}$, where $c_2(\gamma_0)>0$ is an absolute constant. As usual, this implies that a
random $U\in O(n)$ satisfies
\begin{equation}Z_q^{\circ }(\mu )\cap U(Z_q^{\circ }(\mu ))\subseteq \frac{c}{\sqrt{q}}\,B_2^n,\end{equation}
with probability greater than $1-e^{-n}$, where $c>0$ is an absolute constant. This estimate appears in
\cite{Klartag-EMilman-2012b} (and a second proof is given in \cite{Brazitikos-Stavrakakis-2014}).

Assuming that $\sup_nL_n\ls\alpha $ we may apply the same reasoning for every $1\ls q\ls n$:
choosing $k=\gamma_0n$ we see from \eqref{eq:Rqk-variant-3} and \eqref{eq:Rqk-variant-4} that
\begin{equation}c_{\gamma_0n}(Z_q^{\circ }(\mu ))\ls c_1(\gamma_0)\frac{\alpha }{\sqrt{q}},\end{equation}
where $c_1(\gamma_0)>0$ is an absolute constant. Then, we apply \eqref{eq:LPT} with $s=n/2$ and $m=(1-\gamma_0)n$
to get that a random subspace $E\in G_{n,n/2}$ satisfies
\begin{equation}R(Z_q^{\circ }(\mu )\cap E)\ls c_3\cdot c_{\gamma_0n}(Z_q^{\circ }(\mu ))\ls c_2(\gamma_0)\frac{\alpha }{\sqrt{q}}\end{equation}
with probability greater than $1-2e^{-n/4}$, where $c_2(\gamma_0)>0$ is an absolute constant. Finally, this implies that a
random $U\in O(n)$ satisfies
\begin{equation}Z_q^{\circ }(\mu )\cap U(Z_q^{\circ }(\mu ))\subseteq \frac{c\alpha }{\sqrt{q}}\,B_2^n,\end{equation}
with probability greater than $1-e^{-n}$, where $c>0$ is an absolute constant.
\end{remark}

\subsection{Random sections of bodies with maximal isotropic constant}

Starting with an isotropic symmetric convex body $K$ in ${\mathbb R}^n$ we can use the method of this section in order to estimate the quantities
\begin{equation}c_t(K)=\min\{R(K\cap F) : F \in G_{n,n-t}\}\end{equation}
for every $t = 0,\ldots,n-1$. From \eqref{eq:not-9} we have
\begin{equation}|K\cap E|^{\frac{1}{n-t}} \ls
c_2\frac{L_{\overline{K_{k+1}}(\pi_{E^{\perp }}(\mu_{K}))}}{L_K}\ls \frac{c_3L_{n-t}}{L_K}
\end{equation}
for every $E\in G_{n,t}$, therefore
\begin{equation}w_t(K)\ls c_4\sqrt{t}\left (\frac{c_3L_{n-t}}{L_K}\right )^{\frac{n-t}{t}}.\end{equation}
Assume that $K$ has maximal isotropic constant, i.e. $L_K=L_n^{\prime }$ (the same argument works if we assume that $L_K$ is
almost maximal, i.e. $L_K\gr \beta L_n^{\prime }$ for some absolute constant $\beta\in (0,1)$). It is known that
$L_{n-t}\ls c_1L_n\ls c_2L_n^{\prime }$ for all $1\ls t\ls n-1$, where $c_1,c_2>0$ are absolute constants. Therefore, we get:

\begin{lemma}\label{lem:max-Ln-1}Let $K$ be an isotropic symmetric convex body in ${\mathbb R}^n$ such that $L_K=L_n^{\prime }$,
and let $1\ls t\ls \lfloor n/2 \rfloor$. Then,
\begin{equation}c_{2t}(K)\ls c_1^{\frac{n-t}{t}}\frac{n}{\sqrt{t}} \log\Big(e + \frac{n}{t}\Big),\end{equation}
where $c>0$ is an absolute constant.
\end{lemma}

Then, we apply \eqref{eq:LPT} with $s=n/2$ and $m=(1-\gamma_0)n$
to get that a random subspace $E\in G_{n,n/2}$ satisfies
\begin{equation}R(K\cap E)\ls c_3\cdot c_{\gamma_0n}(K)\ls c_1(\gamma_0)\sqrt{n}\end{equation}
with probability greater than $1-2e^{-n/4}$, where $c_1(\gamma_0)>0$ is an absolute constant.

Also, since $c_{\gamma_0n}(K)\ls c(\gamma_0)\sqrt{n}$, we may apply \eqref{eq:Rud-Ver} to get:

\begin{theorem}\label{th:max-Ln-2}Let $K$ be an isotropic symmetric convex body in ${\mathbb R}^n$ with $L_K=L_n^{\prime }$. A random $U\in O(n)$ satisfies
\begin{equation}K\cap U(K)\subseteq (c_3\sqrt{n})\,B_2^n,\end{equation} with probability greater than $1-e^{-n}$, where $c_3>0$ is an absolute constant.
\end{theorem}

We can also prove the local analogue of this fact: random proportional sections of a body with maximal isotropic constant
have bounded isotropic constant.

\begin{theorem}\label{th:max-Ln-3}Let $K$ be an isotropic symmetric convex body in ${\mathbb R}^n$ with $L_K=L_n^{\prime }$. A random $F\in G_{n,n/2}$ satisfies
\begin{equation}L_{K\cap F}\ls c_4\end{equation} with probability greater than $1-e^{-c_5n}$, where $c_4,c_5>0$ are absolute constants.
\end{theorem}

\noindent {\it Proof.} It was proved in \cite{Dafnis-Paouris-2010} (see also \cite[Lemma 6.3.5]{BGVV-book}) that if $L_K=L_n^{\prime }$ then
\begin{equation}|K\cap F|^{\frac{1}{n}}\gr c_6\end{equation}
for every $G_{n,n/2}$, where $c_6>0$ is an absolute constant. Since $R(K\cap F)\ls c_3\sqrt{n}$ for a random $F\in G_{n,n/2}$,
for all these $F$ we get
\begin{equation}\frac{n}{2}L_{K\cap F}^2\ls \frac{1}{|K\cap F|^{1+\frac{2}{n}}}\int_{K\cap F}\|x\|_2^2dx
\ls \frac{1}{|K\cap F|^{\frac{2}{n}}}R^2(K\cap F)\ls c_6^{-2}c_3^2n,\end{equation}
which implies that
\begin{equation}L_{K\cap F}\ls c_4,\end{equation}
where $c_4=\sqrt{2}c_6^{-1}c_3$. \prend

\section{Sub-Gaussian subspaces}

In this section we prove Theorem \ref{th:1.3}. We will use E.~Milman's estimates \cite{EMilman-2014} on the mean width $w(Z_q(K))$ of
the $L_q$-centroid bodies $Z_q(K)$ of an isotropic convex body $K$ in ${\mathbb R}^n$.

\begin{theorem}[E.~Milman]\label{th:Emanuel2}Let $K$ be an isotropic convex body in ${\mathbb R}^n$. Then, for all $q\gr 1$ one has
\begin{equation}w(Z_q(K)) \ls c_1\log (1+q)\max\left\{\frac{q\log (1+q)}{\sqrt{n}},\sqrt{q}\right\}L_K\end{equation}
where $c_1>0$ is an absolute constant.
\end{theorem}

We also use the next fact on the diameter of $k$-dimensional projections of symmetric convex bodies (see \cite[Proposition 5.7.1]{AGA-book}).

\begin{proposition}\label{prop:diam-rdm-proj}
Let $D$ be a symmetric convex body in $\mathbb R^n$ and let $1\ls k<n$ and $\alpha >1$.
Then there exists a subset $\Gamma_{n,k}\subset G_{n,k}$  with measure $\nu_{n,k}(\Gamma_{n,k})\gr 1-e^{-c_2\alpha^2k}$ such that
the orthogonal projection of $D$ onto any subspace $F\in \Gamma_{n,k}$ satisfies
\begin{equation}
R(P_F(D))\ls c_3\alpha \max\{w(D),R(D)\sqrt{k/n}\},
\end{equation} where $c_2>0,c_3>1$ are absolute constants.
\end{proposition}

Combining Proposition \ref{prop:diam-rdm-proj} with Theorem \ref{th:Emanuel2} and the fact that $R(Z_q(K))\ls cqL_K$, we get:

\begin{lemma}\label{lem:main-1}Let $K$ be an isotropic convex body in ${\mathbb R}^n$. Given $1\ls q\ls n$ define $k_0(q)$ by the equation
\begin{equation}k_0(q)=\log^2(1+q)\max\{ \log^2(1+q),n/q\}.\end{equation}
Then, for every $1\ls k\ls k_0(q)$, a random $F\in G_{n,k}$ satisfies
\begin{equation}R(P_F(Z_q(K))) \ls c_1\alpha\log (1+q)\max\left\{\frac{q\log (1+q)}{\sqrt{n}},\sqrt{q}\right\}L_K\end{equation}
with probability greater than $1-e^{-c_2\alpha^2k_0(q)}$, where $c_1,c_2>0$ are absolute constants.
\end{lemma}

\noindent {\it Proof.} Since $R(Z_q(K))\ls cqL_K$ we see that
\begin{align}\frac{R(Z_q(K))\sqrt{k_0(q)}}{\sqrt{n}} &\ls \frac{cq}{\sqrt{n}}\log (1+q)\max\left\{ \log (1+q),\frac{\sqrt{n}}{\sqrt{q}}\right\}L_K\\
\nonumber &= c\log (1+q)\max\left\{\frac{q\log (1+q)}{\sqrt{n}},\sqrt{q}\right\}L_K.\end{align}
From Theorem \ref{th:Emanuel2} we have an upper bound of the same order for $w(Z_q(K))$. Then, we apply Proposition \ref{prop:diam-rdm-proj}
for $Z_q(K)$. \prend

\begin{remark}\label{rem:heredit}\rm Note that if $1\ls s\ls k$ then the conclusion of Proposition \ref{prop:diam-rdm-proj} continues to hold
for a random $F\in G_{n,s}$ with the same probability on $G_{n,s}$; this is an immediate consequence of Fubini's
theorem and of the fact that $R(P_H(D))\ls R(P_F(D))$ for every $s$-dimensional subspace $H$ of a $k$-dimensional
subspace $F$ of ${\mathbb R}^n$.
\end{remark}

\noindent {\bf Proof of Theorem \ref{th:1.3}.} We define $q_0$ by the equation
\begin{equation}q_0\log^2(1+q_0)=n.\end{equation}
Note that $q_0\simeq n/(\log n)^2$ and $\log (1+q_0)\simeq \log n$. For every $2\ls q\ls q_0$ we have $q\log^2(1+q)\ls n$, therefore
\begin{equation}k_0(q)= \frac{n\log^2(1+q)}{q}\gr \frac{c_1n\log^2(1+q_0)}{q_0}\end{equation}
for some absolute constant $c_1>0$, because $q\mapsto \log^2(1+q)/q$ is decreasing for $q\gr 4$. It follows that
\begin{equation}k_0(q)\gr c_1\log^4(1+q_0)\gr c_2(\log n)^4\end{equation}
for all $2\ls q\ls q_0$.

Now, we fix $\alpha >1$ and define
\begin{equation}k_0=c_1\log^4(1+q_0).\end{equation}
Using Lemma \ref{lem:main-1} and Remark \ref{rem:heredit}, for every $q\ls q_0$ we can find a set $\Gamma_q\subseteq G_{n,k_0}$
with $\nu_{n,k_0}(\Gamma_q)\gr 1-e^{-c\alpha^2k_0}$ such that
\begin{equation}R(P_F(Z_q(K))) \ls c_3\alpha\log (1+q)\max\left\{\frac{q\log (1+q)}{\sqrt{n}},\sqrt{q}\right\}L_K\ls c_3\alpha\sqrt{q}\log (1+q)L_K\end{equation}
for all $F\in G_{n,k_0}$. If $\Gamma :=\bigcap_{s=1}^{\lfloor\log_2q_0\rfloor }\Gamma_{2^s}$, then
\begin{equation}\nu_{n,k_0} \big(G_{n,k_0}\setminus \Gamma \big)\ls\nu_{n,k_0} \Big(G_{n,k_0}\setminus \bigcap_{s=1}^{\lfloor\log_2n\rfloor }\Gamma_{2^s}\Big)\ls c(\log n)e^{-c\alpha^2k_0}\ls \frac{1}{n^{\log^3n}}\end{equation}
if $\alpha\simeq 1$ is chosen large enough. Then for every $F\in \Gamma $, for all $\theta\in S_F$ and for every $1\ls s\ls \lfloor\log_2q_0\rfloor $ we have
\begin{equation}\label{eq:small-q}\frac{h_{Z_{2^s}(K)}(\theta )}{\sqrt{2^s}}
=\frac{h_{P_F(Z_{2^s}(K))}(\theta )}{\sqrt{2^s}}\ls c_3\alpha\log (1+2^s)L_K\ls c_4\alpha (\log n)L_K.\end{equation}
Taking into account the fact that if $2^s\ls q<2^{s+1}$ then
\begin{equation}\frac{h_{Z_{q}(K)}(y)}{\sqrt{q}}\ls \frac{h_{Z_{2^{s+1}}(K)}(y)}{2^{s/2}}= \sqrt{2}\frac{h_{Z_{2^{s+1}}(K)}(y)}{2^{(s+1)/2}},\end{equation}
we see that
\begin{equation}\label{eq:final-11}\frac{h_{Z_{q}(K)}(y)}{\sqrt{q}}\ls c_5\alpha (\log n)L_K\end{equation}
for every $F\in \Gamma $, for all $\theta\in S_F$ and for every $2\ls q\ls q_0$.

Next, observe that if $q_0\ls q\ls n$ then we may write
\begin{align}\frac{h_{Z_{q}(K)}(y)}{\sqrt{q}} &\ls \frac{c_6q}{q_0}\frac{h_{Z_{q_0}(K)}(y)}{\sqrt{q}}= \frac{c_6\sqrt{q}}{\sqrt{q_0}}\frac{h_{Z_{q_0}(K)}(y)}{\sqrt{q_0}}
\ls \frac{c_6\sqrt{n}}{\sqrt{q_0}}\frac{h_{Z_{q_0}(K)}(y)}{\sqrt{q_0}}\\
\nonumber &=c_6\log (1+q_0)\frac{h_{Z_{q_0}(K)}(y)}{\sqrt{q_0}}\ls c_7(\log n)\frac{h_{Z_{q_0}(K)}(y)}{\sqrt{q_0}},\end{align}
and hence
\begin{equation}\label{eq:final-2}\frac{h_{Z_{q}(K)}(y)}{\sqrt{q}}\ls c_7\alpha (\log n)^2L_K\end{equation}
for every $F\in \Gamma $, for all $\theta\in S_F$ and for every $q_0\ls q\ls n$.

Recall that $\Psi_2(K)$ is the convex body with support function $h_{\Psi_2(K)}(y)=\|\langle \cdot ,y\rangle\|_{L_{\psi_2}(K)}$. One also has
\begin{equation}h_{\Psi_2(K)}(y)\simeq \sup_{q\gr 2}\frac{h_{Z_q(K)}(y)}{\sqrt{q}}\simeq \sup_{2\ls q\ls n}\frac{h_{Z_q(K)}(y)}{\sqrt{q}}\end{equation}
because $h_{Z_q(K)}(y)\simeq h_{Z_n(K)}(y)$ for all $q\gr n$. Then, \eqref{eq:final-11} and \eqref{eq:final-2} and the fact that $\alpha\simeq 1$ show that
\begin{equation}\|\langle \cdot ,\theta\rangle\|_{L_{\psi_2}(K)}\ls C(\log n)^2L_K\end{equation}
for every $F\in \Gamma $ and  for all $\theta\in S_F$, where $C>0$ is an absolute constant. $\hfill\Box $

\bigskip

\bigskip

\footnotesize
\bibliographystyle{amsplain}

\bigskip

\bigskip

\thanks{\noindent {\bf Keywords:}  Convex bodies, isotropic position, log-concave measures,
centroid bodies, diameter of random sections, sub-Gaussian estimates.}

\smallskip

\thanks{\noindent {\bf 2010 MSC:} Primary 52A23; Secondary 46B06, 52A40, 60D05.}

\bigskip

\bigskip

\noindent \textsc{Apostolos \ Giannopoulos}: Department of
Mathematics, University of Athens, Panepistimioupolis 157-84,
Athens, Greece.

\smallskip

\noindent \textit{E-mail:} \texttt{apgiannop@math.uoa.gr}

\bigskip

\noindent \textsc{Labrini \ Hioni}: Department of
Mathematics, University of Athens, Panepistimioupolis 157-84,
Athens, Greece.

\smallskip

\noindent \textit{E-mail:} \texttt{lamchioni@math.uoa.gr}

\bigskip

\noindent \textsc{Antonis \ Tsolomitis:} Department of Mathematics,
University of the Aegean, Karlovassi 832\,00, Samos, Greece.

\smallskip

\noindent {\it E-mail:} \texttt{atsol@aegean.gr}

\end{document}